\input amstex
\documentstyle{amsppt}
\magnification=1050
\input xypic
\NoBlackBoxes
\NoRunningHeads
\hsize 6.5truein
\document

\define\dfn#1{\definition{\bf\underbar{Definition #1}}}

\define\textit{\it}
\font\script = rsfs10

\def\tbrm#1{\text{\bf{#1}}}

\def\aA{\hbox{\script A}}

\topmatter
\title
The Baum-Connes assembly map and the generalized Bass Conjecture
\endtitle
\vskip.2in
\author
C. Ogle (OSU)
\endauthor
\date
May 2007
\enddate
\email
ogle@math.mps.ohio-state.edu, ronji@math.iupui.edu
\endemail
\endtopmatter
\vskip.5in

\head Introduction
\endhead
\vskip.2in

In the early 1980's, P. Baum and A. Connes defined an assembly map
$$
\aA_*^{G,a}: KK^G_*(C(\underline{E}G),\Bbb C)\to K^t_*(C^*_r(G))
\tag0.1
$$
where $G$ denotes a locally compact group, $\underline{E}G$ the classifying space for proper $G$-actions, $C(\underline{E}G)$ the $G$-algebra of complex-valued functions on $\underline{E}G$ vanishing at infinity, and $KK^G_*(C(\underline{E}G),\Bbb C)$ the $G$-equivariant $KK$-groups of        
$(\underline{E}G)$ with coefficients in $\Bbb C$, while $K_*^t(C^*_r(G))$ represents the topological $K$-groups of the reduced $C^*$-algebra of $G$. The original details of this map appeared (a few years later) in [BC1] and [BC2], with further elaborations in [BCH]. As shown in [BC3], when $G$ is discrete the left-hand side admits a Chern character which may be represented as
$$
ch_*^{BC}(G): KK^G_*(C(\underline{E}G),\Bbb C)\to \underset{x\in\ fin(<G>)}\to{\oplus} H_*(BG_x;\Bbb C)\otimes HPer_*(\Bbb C)
$$
where $fin(<G>)$ is the set of conjugacy classes of $G$ corresponding to elements of finite order, $G_x$ the centralizer of $g$ in $G$ where $x = <g>$, and $HPer_*(\Bbb C)$ the periodic cyclic homology of $\Bbb C$. Note that $H_*(BH;\Bbb C)\otimes HPer_*(\Bbb C)$ are simply the $2$-periodized complex homology groups of $BH$, and (via the classical Atiyah-Hirzebruch Chern character) can be alternatively viewed as the complexified $K$-homology groups of $BH$. Upon complexification, the map $ch_*^{BC}(G)$ is an isomorphism. The original construction of Baum and Connes $\aA_*^{G,a}$ was analytical. Motivated by the need to construct a homotopical analogue to their map, we constructed an assembly map in [O1] which we will denote here as
$$
\aA_*^{G,h}\otimes\Bbb C: H_*(\underset{x\in\ fin(<G>)}\to\coprod BG_x;\tbrm{K}(\Bbb C))\otimes\Bbb C\to K_*^t(C^*_r(G))\otimes\Bbb C
$$
where $\tbrm{K}(\Bbb C)$ denotes the 2-periodic topological $K$-theory spectrum of $\Bbb C$. The construction of this map amounted to an extension of the classical assembly map constructed in [L] which was designed to take into account the contribution coming from the conjugacy classes of finite order. The two essential features of $\aA_*^{G,h}\otimes\Bbb C$, shown in [O1], were (i) it factors through $K^t_*(\Bbb C[G])\otimes\Bbb C$ (where $K^t_*(\Bbb C[G])$ denotes the \underbar{Bott-periodized} topological $K$-theory of the complex group algebra, topologized with the fine topology), and (ii) the composition of $\aA_*^{G,h}\otimes\Bbb C$ with the complexified Chern-Connes-Karoubi-Tillmann character $ch^{CK}_*: K_*(\Bbb C[G])\otimes\Bbb C\to HC_*(\Bbb C[G])$ was effectively computable (see below). What we did not do in [O1] was show that $\aA_*^{G,a}\otimes\Bbb C$ and $\aA_*^{G,h}\otimes\Bbb C$ agree. Since this initial work, there have been numerous extensions and reformulations of the Baum-Connes assembly map, as well as of the original Baum-Connes conjecture, which states that the map in (0.1) is an isomorphism. These extensions typically are included under the umbrella term \lq\lq Isomorphism Conjecture\rq\rq, (formulated for both algebraic and topological $K$-theory; cf. [DL], [FJ], [LR]). Thanks to [HP], we now know that the different formulations of these assembly maps (e.g., homotopy-theoretic vs. analytical) agree.
\vskip.2in

Abbreviating $KK^G_*(C(\underline{E}G),\Bbb C)$ as $K^G_*(\underline{E}G)$ (read: the equivariant $K$-homology of the proper $G$-space $\underline{E}G$), our main result is 

\proclaim{\bf\underbar{Theorem 1}} There is a commuting diagram
$$
\diagram
K^G_*(\underline{E}G)\rto^{\aA_*^{G,DL}}\dto^{ch^{?}_*} & K^t_*(\Bbb C[G])\dto^{ch^{CK}_*}\\
HC_*^{fin}(\Bbb C[G])\rto|<<\tip & HC_*(\Bbb C[G])
\enddiagram
$$
where $\aA_*^{G,DL}$ is the homotopically defined assembly map of [DL], ${}^{fin}HC_*(\Bbb C[G]) :=  \underset{x\in\ fin(<G>)}\to{\oplus} HC_*(\Bbb C[G])_x\cong  \underset{x\in\ fin(<G>)}\to{\oplus} H(BG_x;\Bbb C)\otimes HC_*(\Bbb C)$ is the \underbar{elliptic summand} of $HC_*(\Bbb C[G])$ [JOR], the lower horizontal map is the obvious inclusion, and the Chern character $ch^?_*$ becomes an isomorphism upon complexification for $*\ge 0$.
\endproclaim

Let $\beta$ denote a bounding class, $(G,L)$ a discrete group equipped with a word-length, and $H_{\beta,L}(G)$ the rapid decay algebra associated with this data [JOR]. We write $K^t_*(H_{\beta,L}(G))$
for the Bott-periodic topological $K$-theory of the topological algebra $H_{\beta,L}(G)$. The \underbar{Baum-Connes} assembly map for $H_{\beta,L}(G)$ is defined to be the composition
$$
\aA_*^{G,\beta}: K^G_*(\underline{E}G)\overset{\aA_*^{G,DL}}\to\longrightarrow K^t_*(\Bbb C[G])\to K^t_*(H_{\beta,L}(G))
\tag{BC}
$$
where the second map is induced by the natural inclusion $\Bbb C[G]\hookrightarrow H_{\beta,L}(G)$.
In [JOR], we conjectured that the image of $ch_*:K^t(H_{\beta,L}(G))\to HC_*^t(H_{\beta,L}(G))$ lies in the elliptic summand ${}^{fin}HC_*^t(H_{\beta,L}(G))$ (conjecture $\beta$-SrBC). As the inclusion $\Bbb C[G]\hookrightarrow H_{\beta,L}(G)$ sends ${}^{fin}HC_*(\Bbb C[G])$ to ${}^{fin}HC_*^t(H_{\beta,L}(G))$, naturality of the Chern character $ch_*^{CK}$ and Theorem 1 implies

\proclaim{\bf\underbar{Corollary 2}} If $\aA_*^{G,\beta}$ is rationally surjective, then $\beta$-SrBC is true.
\endproclaim

Since going down and then across is rationally injective, we also have (compare [O1])

\proclaim{\bf\underbar{Corollary 3}} The assembly map $\aA_*^{G,DL}\otimes\Bbb Q$ is injective for all discrete groups $G$.
\endproclaim

We do not claim any great originality in this paper. In fact, Theorem 1, although not officially appearing in print before this time,  has been a \lq\lq folk-theorem\rq\rq\  known to experts for many years. The connection between the Baum-Connes Conjecture (more precisely a then-hypothetical Baum-Connes-type Conjecture for $\Bbb C[G]$) and the stronger Bass Conjecture for $\Bbb C[G]$ discussed in [JOR] was noted by the author in [O2]. 

There is some overlap of this paper with the results presented in [Ji]. A special case of Theorem 1 (for $* = 0$ and $\Bbb C[G]$ replaced by the $\ell^1$-algebra $\ell^1(G)$) appeared as the main result of [BCM].\vskip.3in

\head
Proof of Theorem 1
\endhead

We use the notation $F_*^{fin}(\Bbb C[G])$ to denote the \underbar{elliptic summand} $\underset{x\in\ fin(<G>)}\to{\oplus}F_*(\Bbb C[G])_x$ of $F_*(\Bbb C[G])$ where $F_*(_-) = HH_*(_-), HN_*(_-), HC_*(_-)$ or $HPer(_-)$.
To maximize consistency with [LR], we write {\bf{S}} for the (unreduced) suspension spectrum of the zero-sphere $S^0$, {\bf{HN}}($R$) resp. {\bf{HH}}($R$) the Eilenberg-MacLane spectrum whose homotopy groups are the \underbar{negative cyclic} resp. \underbar{Hochschild} homology groups of the discrete ring $R$, and $\text{\bf{K}}^a(R)$ the non-connective algebraic $K$-theory spectrum of $R$, with $K_*^a(R)$ representing its homotopy groups. By [LR, diag. 1.6] there is a commuting diagram
$$
\diagram
H_*^G(\underline{E}G;{\text{\bf{S}}})\dto\rrto & & K^a_*(\Bbb Z[G])\dto^{NTr_*}\\
H_*^G(\underline{E}G;{\text{\bf{HN}}}(\Bbb Z))\dto\rto^{\cong} 
& HN_*^{fin}(\Bbb Z[G])\dto\rto|<<\tip & HN_*(\Bbb Z[G])\dto^{h_*} \\
H_*^G(\underline{E}G;{\text{\bf{HH}}}(\Bbb Z))\rto^{\cong} 
& HH_*^{fin}(\Bbb Z[G])\rto|<<\tip& HH_*(\Bbb Z[G])
\enddiagram
\tag1.1
$$
where the top horizontal map is the composition 
$$
H_*^G(\underline{E}G;{\text{\bf{S}}})\to H_*^G(\underline{E}G;{\text{\bf{K}}^a}(\Bbb Z))\overset{\aA^{G,DL}}\to\longrightarrow K_*(\Bbb Z[G])
$$
referred to as the the \underbar{restricted} assembly map for the algegraic $K$-groups of $\Bbb Z[G]$. The other two horizontal maps are the assembly maps for negative cyclic and Hochschild homology respectively. The upper left-hand map is induced by the map from the sphere spectrum to the Eilenberg-MacLane spectrum \text{\bf{HN}}, which may be expressed as the composition of spectra $\text{\bf{S}}\to \text{\bf{K}}^a(\Bbb Z)\to \text{\bf{HN}}$. By [LR], the composition on the left is a rational equivalence.

Let $\Bbb C^{\delta}$ denote the complex numbers $\Bbb C$ equipped with the discrete topology. Tensoring with $\Bbb C$ and combined with the inclusion of group algebras $\Bbb Z[G]\hookrightarrow \Bbb C^{\delta}[G]$, (1.1) yields the commuting diagram
$$
\diagram
H_*^G(\underline{E}G;\Bbb Q)\otimes\Bbb C\dto^{\cong}\rto & K^a_*(\Bbb C^{\delta}[G])\otimes\Bbb C\dto^{NTr_*}\\
HN_*^{fin}(\Bbb C[G])\rto|<<\tip & HN_*(\Bbb C[G])
\enddiagram
\tag1.2
$$
Next, we consider the transformation from algebraic to topological $K$-theory, induced by the map of group algebras $\Bbb C^{\delta}[G]\to \Bbb C[G]$ which is the identity on elements. By the results of [CK], [W] and [T], there is a commuting diagram
$$
\diagram
K_*^a(\Bbb C^{\delta}[G])\otimes\Bbb C\dto^{NTr_*}\rto & K^t_*(\Bbb C[G])\otimes\Bbb C\dto^{ch_*(\Bbb C[G])}\\
HN_*(\Bbb C[G])\rto & HPer_*(\Bbb C[G])
\enddiagram
\tag1.3
$$
where $ch_*(\Bbb C[G])$ is the Connes-Karoubi Chern character for the fine topological algebra $\Bbb C[G]$, and the bottom map is the transformation from negative cyclic to periodic cyclic homology. 

We can now consider our main diagram
$$

\diagram
H_*^G(\underline{E}G;\Bbb C)\otimes K_*(\Bbb C)\rto\dto & K_*^{a}(\Bbb C^{\delta}[G])\otimes\Bbb C\otimes K_*(\Bbb C)\rto\dto & K_*^{t}(\Bbb C[G])\otimes\Bbb C\otimes K_*(\Bbb C)\rto\dto^{ch_*(\Bbb C[G])\otimes ch_*(\Bbb C[\{id\}])} & K^t_*(\Bbb C[G])\otimes\Bbb C\dto^{ch_*(\Bbb C[G])}\\
HN_*^{fin}(\Bbb C[G])\otimes K_*(\Bbb C)\dto\rto & HN_*(\Bbb C[G])\otimes K_*(\Bbb C)\rto & HPer_*(\Bbb C[G])\otimes HPer_*(\Bbb C)\rto^(.6){\cong} & HPer_*(\Bbb C[G])\\
HPer_*^{fin}(\Bbb C[G])\otimes K_*(\Bbb C)\rrto^{\cong} & & HPer_*^{fin}(\Bbb C[G])\otimes HPer_*(\Bbb C)\rto^(.6){\cong}\uto|<<\tip & HPer_*^{fin}(\Bbb C[G])\uto|<<\tip
\enddiagram
\tag1.4
$$

The top left square commutes by (1.2), and the middle top square commutes by (1.3). The upper right square commutes by virtue of the fact that the Connes-Karoubi-Chern character is a homomorphism of graded modules, which maps the $K_*^t(\Bbb C)$-module $K_*^t(\Bbb C[G])$ to the $HPer_*(\Bbb C)$-module $HPer_*(\Bbb C[G])$, with the map of base rings induced by isomorphism $ch_*(\Bbb C[\{id\}]): K_*^t(\Bbb C)\otimes\Bbb C\overset\cong\to\longrightarrow HPer_*(\Bbb C)$. The lower left square commutes trivially, while the lower right commutes by the naturality of the inclusion $HPer_*^{fin}(\Bbb C[G])\hookrightarrow HPer_*(\Bbb C[G])$ with respect to the module structure over $HPer_*(\Bbb C)$.
Summarizing, we get a commuting diagram
$$
\diagram
H_*^G(\underline{E}G;\Bbb C)\otimes K_*(\Bbb C)\rto\dto^{\cong} &  K^t_*(\Bbb C[G])\otimes\Bbb C\dto^{ch_*(\Bbb C[G])}\\
HPer_*^{fin}(\Bbb C[G])\dto\rto|<<\tip & HPer_*(\Bbb C[G])\dto\\
HC_*^{fin}(\Bbb C[G])\rto|<<\tip & HC_*(\Bbb C[G])
\enddiagram
\tag1.5
$$
where the bottom square is induced by the transformation $HPer_*(_-)\to HC_*(_-)$, which respects the summand decomposition indexed on conjugacy classes. Restricted the elliptic summand yields the map$HPer_*^{fin}(\Bbb C[G])\to HC_*^{fin}(\Bbb C[G])$ which is an isomorphism for $*\ge 0$, implying the result stated in Theorem 1.
\vskip.3in

.

\enddocument